\documentstyle[12pt]{article}
\begin{document}
\newcommand{\Si}{\Sigma}
\newcommand{\tr}{{\rm tr}}
\newcommand{\ad}{{\rm ad}}
\newcommand{\Ad}{{\rm Ad}}
\newcommand{\ti}[1]{\tilde{#1}}
\newcommand{\om}{\omega}
\newcommand{\Om}{\Omega}
\newcommand{\de}{\delta}
\newcommand{\al}{\alpha}
\newcommand{\te}{\theta}
\newcommand{\vth}{\vartheta}
\newcommand{\be}{\beta}
\newcommand{\la}{\lambda}
\newcommand{\La}{\Lambda}
\newcommand{\D}{\Delta}
\newcommand{\ve}{\varepsilon}
\newcommand{\ep}{\epsilon}
\newcommand{\vf}{\varphi}
\newcommand{\G}{\Gamma}
\newcommand{\ka}{\kappa}
\newcommand{\ip}{\hat{\upsilon}}
\newcommand{\Ip}{\hat{\Upsilon}}
\newcommand{\ga}{\gamma}
\newcommand{\ze}{\zeta}
\def\bfa{{\bf a}}
\def\bfb{{\bf b}}
\def\bfc{{\bf c}}
\def\bfd{{\bf d}}
\def\bfm{{\bf m}}
\def\bfn{{\bf n}}
\def\bfp{{\bf p}}
\def\bfu{{\bf u}}
\def\bfv{{\bf v}}
\def\bft{{\bf t}}
\newcommand{\li}{\lim_{n\rightarrow \infty}}
\newcommand{\mat}[4]{\left(\begin{array}{cc}{#1}&{#2}\\{#3}&{#4}
\end{array}\right)}
\newcommand{\si}{\sigma}
\newcommand{\beq}[1]{\begin{equation}\label{#1}}
\newcommand{\eq}{\end{equation}}
\newcommand{\beqn}[1]{\begin{eqnarray}\label{#1}}
\newcommand{\eqn}{\end{eqnarray}}
\newcommand{\p}{\partial}
\newcommand{\di}{{\rm diag}}
\newcommand{\oh}{\frac{1}{2}}
\newcommand{\su}{{\bf su_2}}
\newcommand{\uo}{{\bf u_1}}
\newcommand{\GL}{{\rm GL}(N,{\bf C})}
\newcommand{\SL}{{\rm SL}(N,{\bf C})}
\newcommand{\gl}{gl(N,{\bf C})}
\newcommand{\PSL}{{\rm PSL}_2({\bf Z})}
\def\f1#1{\frac{1}{#1}}
\newcommand{\rar}{\rightarrow}
\newcommand{\upar}{\uparrow}
\newcommand{\sm}{\setminus}
\newcommand{\ms}{\mapsto}
\newcommand{\bp}{\bar{\partial}}
\newcommand{\bz}{\bar{z}}
\newcommand{\bA}{\bar{A}}
\newcommand{\sect}[1]{\setcounter{equation}{0}\section{#1}}
\renewcommand{\theequation}{\thesection.\arabic{equation}}
\newtheorem{predl}{Proposition}[section]
\newtheorem{defi}{Definition}[section]
\newtheorem{rem}{Remark}[section]
\newtheorem{cor}{Corollary}[section]
\newtheorem{lem}{Lemma}[section]
\newtheorem{theor}{Theorem}[section]

\vspace{0.3in}
\begin{flushright}
 ITEP-TH29/97\\
 \end{flushright}
\vspace{10mm}
\begin{center}
{\Large\bf Painlev\'{e} - Calogero correpondence.}
\footnote{A talk given by M.O. on a workshop {\em "Calogero-Moser-Sutherland
Models"}, March 97, CRM, Montreal}
\\
\vspace{5mm}
A.M.Levin\\
{\sf Institut of Oceanology, Moscow, Russia,} \\
{\em e-mail andrl@landau.ac.ru}\\

M.A.Olshanetsky
\\
{\sf Institut of Theoretical and Experimental Physics, Moscow, Russia,} \\
{\em e-mail olshanet@heron.itep.ru}\\

\vspace{5mm}
\end{center}
\begin{abstract}
It is proved that the Painlev\'{e} VI equation
$(PVI_{\al,\be,\ga,\de})$  for the special
values of constants $(\al=\frac{\nu^2}{4},\be=-\frac{\nu^2}{4},
\ga=\frac{\nu^2}{4},\de=\f1{2}-\frac{\nu^2}{4})$ is a
reduced hamiltonian
system. Its phase space  is the set of flat $SL(2,{\bf C})$
connections over
elliptic curves with a marked point and time of the system is
given by the elliptic module.
This equation can be derived via reduction procedure from
the free infinite
hamiltonian system. The phase space of later
is the affine space  of smooth connections and the "times are
the Beltrami differentials. This approach allows to define
the associate linear
problem, whose  isomonodromic deformations is provided by the
Painlev\'{e} equation and the Lax pair. 
In addition, it leads to description of solutions by a linear
procedure.
This scheme can be generalized to $G$ bundles over
Riemann curves with marked points, where $G$ is a
simple complex Lie group. In some special limit
such hamiltonian systems convert into the Hitchin systems.
In particular, for $\SL$ bundles
over  elliptic curves with a marked point we obtain in this
limit the elliptic
Calogero $N$-body system. Relations to the classical limit of the
Knizhnik- Zamolodchikov-Bernard equations is discussed.
\end{abstract}

\section {Introduction}
\setcounter{equation}{0}

{\bf 1}.We learned from Yu. Manin's lectures in MPI (Bonn, 1996) about elliptic
form of the famous
Painlev\'{e} VI equation (PVI) \cite{Ma}. In this representation
PVI looks very
similar to the elliptic Calogero-Inozemtsev-Treibich-Verdier (CITV)
 rank one system \cite{Ca,In,TV}. Namely, the both equations
 are hamiltonian with
the same symplectic structure for two dynamical variables and
the same Hamiltonians.
The only difference is that the time in the PVI system is
nothing else as the
elliptic module. Therefore, it is non autonomous hamiltonian system,
while
the CITV Hamiltonian is independent of time. This similarity
is not accidental and based on very closed geometric origin
of the both systems,
which we will elucidate in this talk.
It should be confessed from the very beginning that at the present
time our approach is  cover only
the one parametric family of $PVI_{\al,\be,\ga,\de}$.
This family corresponds to the standard two-body
elliptic Calogero equation.
\bigskip

{\bf 2. Painlev\'{e} VI and Calogero equations.}
The Painlev\'{e} VI $PVI_{\al,\be,\ga,\de}$
equation depends of four free parameters and has the
form
$$
\frac{d^2X}{dt^2}=\frac{1}{2}(\frac{1}{X}+\frac{1}{X-1}+\frac{1}{X-t})
(\frac{dX}{dt})^2-(\frac{1}{t}+\frac{1}{t-1}+
\frac{1}{X-t})\frac{dX}{dt}+
$$
\beq{I.1}
+\frac{X(X-1(X-t)}{t^2(t-1)^2}(\al+\be\frac{t}{X^2}+
\ga\frac{t-1}{(X-1)^2}
+\de\frac{t(t-1)}{X-t)^2})
\eq
It is a hamiltonian systems \cite{O1}, but we will write
the symplectic form and the Hamiltonian below in another variables.
Among some distinguish features of this equation we are interesting in
its relation to the isomonodromic deformations of
linear differential equations.
This approach was investigated  by Fuchs \cite{F}, while
first  $PVI_{\al,\be,\ga,\de}$ was written down by Gambier \cite{G}.
The equation has a lot
of different applications (see \cite{PT}).
We shortly present  $PVI_{\al,\be,\ga,\de}$ in terms
of elliptic functions \cite{Ma}.

Let $\wp(u|\tau)$ be the Weiershtrass function on the elliptic curve
$T^2_{\tau}={\bf C}/({\bf Z}+{\bf Z}\tau)$, and\\
$e_i=\wp(\frac{T_i}{2}|\tau),~(i=1,2,3)~~
(T_0,\ldots,T_3)=(0,1,\tau,1+\tau).$
Consider instead of $(t,X)$ in (\ref{I.1})  the new variables
\beq{I.2}
(\tau,u)\rightarrow (t=\frac{e_3-e_1}{e_2-e_1},X=
\frac{\wp(u|\tau)-e_1}{e_2-e_1}).
\eq
Then $PVI_{\al,\be,\ga,\de}$ takes the form.
\beq{I.3}
\frac{d^2u}{d\tau^2}=\p_uU(u|\tau),~~
U(u|\tau)=\frac{1}{(2\pi i)^2}\sum_{j=0}^3\al_j
\wp(u+\frac{T_j}{2}|\tau),
\eq
$(\al_0,\ldots,\al_3)=(\al,-\be,\ga,\f1{2}-\de)$.
 As usual in non autonomous case, the equations of motion
(\ref{I.3}) are derived from the variations of the degenerated
symplectic form
\beq{I.4}
\om=\de v\de u-\de H\de\tau,~~H=\frac{v^2}{2}+U(u|\tau),
\eq
which is defined over the extended phase space ${\cal P}=\{v,u,\tau\}$.
The semidirect product of ${\bf Z}+{\bf Z}\tau$ and the modular
group act on the dynamical variables $(v,u,\tau)$ 
preserving (\ref{I.4}).

Let us introduce the new parameter $\ka$ (the level) and instead
of (\ref{I.4})
consider
\beq{I.6}
\om=\de v\de u-\f1{\ka}\de H\de\tau.
\eq
It corresponds to the overall rescaling of constants
$\al_j\rar\frac{\al_j}{\ka^2}$.
Put $\tau=\tau_0+\ka t^H$ and consider the system in the limit
$\ka\rar 0$,
which is called the critical level.
We come to the equation
\beq{I.7}
\frac{d^2u}{(dt^H)^2}=\p_uU(u|\tau_0),
\eq
It is just the rank one  $CITV_{\al,\be,\ga,\de}$ equation.
 Thus, we have in this limit\\
$PVI_{\al,\be,\ga,\de}\stackrel{\ka\rar 0}\longrightarrow
CI_{\al,\be,\ga,\de}.$

Consider   one-parametric family $PVI_{\frac{\nu^2}{4},-\frac{\nu^2}{4},
\frac{\nu^2}{4},\f1{2}-\frac{\nu^2}{4}}$
The potential (\ref{I.3}) takes the form
\beq{I.8}
U(u|\tau)=\frac{1}{(4\pi i)^2}\nu^2\wp(2u|\tau),
\eq
We will prove that (\ref{I.6}) with the potential
(\ref{I.8}) describe the dynamic of flat connections
of ${\rm SL}(2,{\bf C})$  bundles over elliptic curves $T_{\tau}$
with one marked
point $\Si_{1,1}$. In fact, $u$ lies on the Jacobian of
$T_{\tau}$, $(v,u)$
 defines a flat bundle, and $\tau$ defines a point in
the moduli space ${\cal M}_{1,1}=\{\Si_{1,1}\}$.
The choice of  the polarization of connections, in other words
$v$ and $u$, depends on
complex structure of $\Si_{1,1}$. The extended phase space ${\cal P}$
includes beside the dynamical variables $v$ and $u$ the "time" $\tau$.
 It is the bundle  over ${\cal M}_{1,1}$ with the fibers
 ${\cal R}=\{v,u\}$,
 which is endowed by
the degenerated symplectic structure
$\om$ (\ref{I.6}).
This system is derived by a reduction
procedure from some free, but infinite hamiltonian system.
In this way we obtain the Lax equations, the linear system which
monodromies preserve by (\ref{I.7}) with $U(u|\tau)$ (\ref{I.8}) and the
explicit solutions of the Cauchy problem via the so-called
projection method. The discrete symmetries of
(\ref{I.4}) are nothing else as the remnant gauge symmetries.
On the critical level it is
just two-body elliptic Calogero system. The corresponding quantum
system is identified with the KZB equation \cite{KZ,B}
for the one-vertex correlator on $T_\tau$.
In the similar way $PVI_{\frac{\nu^2}{4},-\frac{\nu^2}{4},
\frac{\nu^2}{4},\f1{2}-\frac{\nu^2}{4}}$ is the classical limit of the
KZB for $\ka\neq 0$.
\bigskip

{\bf 3.}This particular example has far-reaching generalizations.
Consider a  phase space, which  is
the moduli space of flat connections
${\cal A}$ of $G$ bundle
over Riemann curve $\Si_{g,n}$ of genus $g$ with n marked points, where
$G$ is a complex simple Lie group.
While the flatness is the topological property of bundles,
 the polarization of connections ${\cal A}=(A,\bA)$ depends
on the choice of complex
structure on $\Si_{g,n}$. Therefore, we consider a  bundle ${\cal P}$
 over the moduli space ${\cal M}_{g,n}$ of curves with
flat connections ${\cal R}=(A,\bA)$
as fibers. The fibers are supplemented by elements of coadjoint
orbits ${\cal O}_a$ in the marked points $x_a$.
There exists a closed degenerate two-form $\om$
on ${\cal P}$, which is non degenerate on the fibers.
The equations of motions are defined as variations of the dynamical
variables along the null-leaves of this symplectic form.
We call them as the {\sl hierarchies of isomonodromic deformations}
(HID). They are attended by the {\sl Whitham hierarchies},
which has occurred earlier
in \cite{Kr1} as a result of the averaging procedure, and then in
\cite{DVV} as the classical limit of "string equations".
Our approach is closed to the Hitchin construction of integrable systems,
living on the cotangent
bundles to the moduli space of holomorphic $G$ bundles \cite{H1},
generalized for singular curves in \cite{Ne,ER}.
Namely, the connection $\bA$ plays the same role as in the
Hitchin scheme, while $A$ replaces the Higgs field.
Essentially, our construction is local - we
work over a vicinity of some fixed curve $\Si_{g,n}$ in
${\cal M}_{g,n}$.
The coordinates of tangent vector to ${\cal M}_{g,n}$ in this point
play role of times, while the Hitchin times have nothing to do
with the moduli space. The Hamiltonians are the same quadratic Hitchin
Hamiltonians, but now they are time dependent. There is
a free parameter  $\ka$ ({\sl the level}) in our construction.
On the critical level  $(\ka=0)$, after rescaling the times, our systems
convert into the Hitchin systems. In concrete examples our
work is based essentially on \cite{Ne}, which deals with
the same systems on the critical level.

 As the Hitchin systems, HID can be derived by the
symplectic reduction from a free infinite  hamiltonian system.
In our case
the upstairs extended phase space is the
space of the affine connections and
the Beltrami differentials. We consider its symplectic
quotient with respect of gauge action on the connections. 
In addition, to come to the moduli space ${\cal M}_{g,n}$
we need the
subsequent factorization under the action of the diffeomorphisms,
which effectively acts on the Beltrami differentials only.
Apart from
the last step, this derivation resembles the construction
of the KZB systems
in \cite{ADPW}, where they are derived as a quantization
of the very similar symplectic quotient.
Our approach allows to write down
 the Lax pairs, prove that the HID are consistency conditions of the
isomonodromic deformations of the
linear Lax equations, and, therefore, justify the notion HID.
Moreover, we describe solutions via linear procedures (the projection
method). HID are the quasi classical limit of the KZB equations for
$(\ka\neq 0$, as the Hitchin systems are the quasi classical
limit of the KZB
equations
on the critical level \cite{Ne,I}. The quantum counterpart
of the Whitham
hierarchy is the flatness condition, which discussed in \cite{H2} within
derivation of KZB. The interrelations between quantizations of
isomonodromic
deformations and the KZB eqs were discussed in \cite{R,Ha,Ko}.

For genus zero our procedure leads to Schlesinger equations.
We restrict ourselves to simplest cases
with only simple poles of connections. Therefore, we don't
include in the phase space the Stokes parameters.
This phenomen was investigated in the rational case in detail
in \cite{JMU}. For
genus one we obtain a particular case of the Painlev\'{e} VI equation (for
$SL(2,{\bf C}$ bundles with one marked point),
 generalization of this case  on arbitrary simple groups and arbitrary
number of marked points.
\bigskip

{\bf Acknowledgments.}\\
{\sl  We are
thankful to Yu.Manin - his lectures and discussions with
him concerning PVI, stimulated
our interests to these problems. 
 We are grateful to the Max-Planck-Institut f\"{u}r Mathemamatik in Bonn
for the hospitality, where this work was started. 
We would like to thank our collegues V.Fock,  A.Losev,
A.Morozov, N.Nekrasov, and A.Rosly for fruitfull discussions duering
working on this subject.
 The work is
supported in part by by Award No.
RM1-265 of the US Civilian Research \& Development Foundation
(CRDF) for the Independent States of the Former Soviet Union,
and grant  96-15-96455  for support of scientific schools  (A.L);
grants RFBR-96-02-18046, Award No.
RM2-150 of the US Civilian Research \& Development Foundation
(CRDF) for the Independent States of the Former Soviet Union, INTAS
930166 extension,
 and grant  96-15-96455  for support of scientific schools (M.O).}

\section{Symplectic reduction}
\setcounter{equation}{0}

{\bf 1. Upstairs extended phase space.}
 Let $\Si_{g,n}$ be
a Riemann curve  of genus $g$ with
$n$ marked points. Let us fix the complex structure of $\Si_{g,n}$
defining
local coordinates $(z,\bz)$ in open maps covering $\Si_{g,n}$.
Assume that the marked points $(x_1,\ldots,x_n)$
are in the generic positions. The deformations
of the basic complex structure are determined by the
Beltrami differentials $\mu$,
which are smooth $(-1,1)$ differentials on $\Si_{g,n}$, 
$\mu\in{\cal A}^{(-1,1)}(\Si_{g,n})$ .
We identify this set with the space of times ${\cal N}'$.
 The Beltrami  differentials  can be defined in the following way.
Consider the  chiral diffeomorphisms of $\Si_{g,n}$
\beq{2a.2}
w=z-\ep(z,\bz),~~\bar{w}=\bz
\eq
and the corresponding one-form $dw$.  Up to the conformal factor
$1-\p \ep(z,\bz)$, it is equal
\beq{1a.2}
dw=dz-\mu d\bz,~~\mu=\frac{\bp \ep(z,\bz)}{1-\p \ep(z,\bz)}.
\eq
The new holomorphic structure is defined by the deformed
antiholomorphic operator annihilating $dw$, while the
antiholomorphic structure is kept unchanged
$$\p_{\bar{w}}=\bp+\mu\p,~~\p_w=\p.$$
In addition, assume that $\mu$ vanishes
in the marked points $\mu(z,\bz)|_{x_a}=0.$
We consider small deformations of the basic complex
structure $(z,\bz)$. It allows  to replace (\ref{1a.2}) by
\beq{2b.2}
\mu=\bp \ep(z,\bz).
\eq

 Let ${\cal E}$ be a principle stable
$G$ bundle over a Riemann curve $\Si_{g,n}$. Assume that $G$
is a complex simple Lie group.
The phase space ${\cal R}'$ is recruited
from the following data:\\
i)the affine space $\{{\cal A}\}$ of Lie$(G)$-valued
connection on ${\cal E}$.\\
It has the following component description:\\
a)   $C^{\infty}$ connection $\{\bar{A}\}$ ,
corresponding to the $d\bar{w}=d\bz$ component of ${\cal A}$;\\
b)The dual  to the previous space the space $\{A\}$ of $dw$
components of ${\cal A}$.
 $A$ can have simple poles in the marked points. Moreover,
assume that $\bA\mu$ is a  $C^{\infty}$ object ;\\
ii)cotangent bundles
 $T^*G_a=\{(p_a,g_a),~p_a\in {\rm Lie}^*(G_a),~g_a\in G_a\},
 ~(a=1,\ldots,n)$ in the points
$(x_1,\ldots,x_n)$.\\
There is the canonical  symplectic form on ${\cal R}'$
\beq{1.2}
\om_0=\int_{\Si}<\de A,\de\bar{A}>+
2\pi i\sum_{a=1}^n\de<p_a,g_a^{-1}\de g_a>,
\eq
where $<~,~>$ denotes the Killing form on Lie$(G)$.
\bigskip

Consider the bundle ${\cal P}'$ over ${\cal N}'$ with
${\cal R}'$ as the fibers. It plays role of the extended phase space.
There exists the degenerate form  on  ${\cal P}'$
\beq{2.2}
\om=\om_0-
\frac{1}{\ka}
\int_{\Si}<\de A, A>\de \mu.
\eq
Thus, we deal with the infinite set of Hamiltonians
$< A, A>(z,\bz)$,
parametrized by points of $\Si_{g,n}$ and corresponding set of
times $\mu(z,\bz)$.
The equations of motion.
 take the form
\beq{3.2}
\frac{\p A}{\p\mu}(z,\bar{z}) =0,~~
\ka \frac{\p \bA}{\p\mu}(z,\bar{z})=A(z,\bar{z}),~~
\frac{\p p_b}{\p\mu}=0,~~
\frac{\p g_b}{\p \mu}=0.
\eq
We will apply the formalism of hamiltonian reduction to
these systems.

\bigskip

 {\bf 2. Symmetries}.
The form $\om$ (\ref{2.2}) is invariant with respect to
the action of
the group ${\cal G}_0$ of  diffeomorphisms of $\Si_{g,n}$,
which are trivial in vicinities ${\cal U}_a$ of marked points:
\beq{7.2}
{\cal G}_0=\{z\rar N(z,\bar{z}),\bar{z}\rar \bar{N}(z,\bar{z}),),~
 N(z,\bar{z})=z+o(|z-x_a|),~z\in {\cal U}_a\}.
\eq

Another infinite gauge symmetry of  the form (\ref{2.2})   is the
group
${\cal G}_1=\{f(z,\bar{z})\in C^\infty (\Si_{g},G)\}$
that acts on the dynamical fields as
$$
A\rar f(A+\ka\p)f^{-1},~~~\bar{A}\rar f(\bar{A}+\bp+\mu\p)f^{-1},
$$
\beq{8.2}
(\bA'\rar f(\bar{A}'+\bp)f^{-1}),
\eq
$$
p_a\rar f_ap_af^{-1}_a,
~~g_a\rar g_af^{-1}_a,~~(f_a=\lim_{z\rar x_a}f(z,\bz)),
~~
\mu\rar\mu.
$$
In other words, the gauge action of ${\cal G}_1$
 does not touch the base  ${\cal N}'$ and transforms only the
fibers ${\cal R}'$.
The whole gauge group is the semidirect product
${\cal G}_1\oslash{\cal G}_0.$
\bigskip

{\bf 3. Symplectic reduction with respect to ${\cal G}_1$.}
Since the symplectic form (\ref{2.2}) is closed
(though is degenerated)
one can consider the symplectic quotient
of the extended phase space ${\cal P}'$
under the action of the gauge transformations   (\ref{8.2}).
They  are generated by the moment constraints
\beq{10.2}
F_{A,\bA}(z,\bz)-2\pi i\sum_{a=1}^n\de^2(x_a)p_a=0,
\eq
where
$F_{A,\bA}=(\bp +\p\mu)A-\ka\p\bA+[\bA.A].$
It means that we deal with the flat connection everywhere on
$\Si_{g,n}$ except
the marked points. The holonomies of $(A,\bA)$ around
the marked points are conjugate to  $\exp 2\pi ip_a$.

Let $(L,{\bar L})$ be the gauge transformed connections
\beq{11.2}
\bA=f(\bar{L}+\bp+\mu\p)f^{-1},~~A=f(L+\ka\p)f^{-1},
\eq
Then (\ref{10.2}) takes the form
\beq{12a.2}
(\bp+\p\mu)L-\ka\p {\bar L}+[\bar{L},L]=2\pi i\sum_{a=1}^n\de^2(x_a)p_a.
\eq
\begin{rem}
The gauge fixing allows to choose $\bA$ in a such way that
$\p\bar L=0$.
Then (\ref{12a.2}) takes the form
 \beq{14.2}
(\bp+\p\mu)L+[\bar{L},L]=\sum_{a=1}^n\de^2(x_a)p_a.
\eq
It coincides with the moment equation for the Hitchin systems
on singular curves
\cite{Ne}.
\end{rem}
We can rewrite (\ref{14.2}) as
\beq{13.2}
\p_{\bar w}L+[\bar{L},L]=
2\pi i\sum_{a=1}^n\de^2(x_a)p_a.
\eq
Anyway, by choosing $\bar{L}$ we fix somehow the gauge
in generic case.
There is additional gauge freedom $h_a$ in the points $x_a$, which
  acts on $T^*G_a$ as \mbox{$p_a\rar p_a$, $g_a\rar h_ag_a$}.
It allows to fix $p_a$ on some coadjoint
orbit $p_a=g_a^{-1}p_a^{(0)}g_a$ and obtain the symplectic quotient
${\cal O}_a=T^*G_a//G_a$. Thus in (\ref{12a.2}) or (\ref{14.2})
$p_a$ are
elements of ${\cal O}_a$.

Let ${\cal I}_{g,n}$ be the equivalence classes of the connections
$(A,\bA)$ with respect to the gauge action (\ref{11.2}) -
 the moduli space of stable flat $G$ bundles
over $\Si_{g,n}$ . It is a smooth
finite dimensional space.
Fixing the conjugacy classes of holonomies $(L,\bar L)$
  around marked points (\ref{12a.2})
 amounts to choose a symplectic leave ${\cal R}$ in ${\cal I}_{g,n}$.
Thereby we come to the symplectic quotient
$$
{\cal R}={\cal R}'//{\cal G}_1=
{\cal J}^{-1}_1(0)/{\cal G}_1\subset{\cal I}_{g,n}.
$$
 The connections $(L,\bar{L})$ in addition to
 $\bfp=(p_1,\ldots,p_n)$ depend
 on  a finite even number of free parameters $2r$ $
(\bfv,\bfu),~\bfv=(v_1,\ldots,v_{r}),~{\bf u}=(u_1,\ldots,u_{r}).$
$$
r=
\left\{
\begin{array}{ll}
0&~g=0,\\
{\rm rank} G,&~g=1\\
(g-1)\dim G,&~g\geq 2.\\
\end{array}
\right.
$$
The fibers ${\cal R}$ are symplectic manifolds with the nondegenerate
symplectic form which is the reduction of (\ref{2.2})
\beq{15.2}
\om_0=\int_{\Si}<\de L,\de \bar{L}>+
2\pi i\sum_{a=1}^n\de<p_a,g_a\de g_a^{-1}>.
\eq

On this stage we come to the bundle ${\cal P}''$ with
the finite-dimensional fibers\\
 ${\cal R}$ over the infinite-dimensional
base ${\cal N}'$  with the symplectic form
\beq{16.2}
\om=\om_0-
\f1{\ka}\int_{\Si}<L,\de L>\de\mu.
\eq
\bigskip

{\bf 4. Factorization with respect to the diffeomorphisms}
${\cal G}_0$.
We can utilize invariance of $\om$ with respect to ${\cal G}_0$
and reduce ${\cal N}'$ to the finite-dimensional
space ${\cal N}$, which is isomorphic to the moduli space
${\cal M}_{g,n}$.
The crucial point is that for the flat connections
the action of diffeomorphisms ${\cal G}_0$ on the connection fields
is generated by the gauge transforms ${\cal G}_1$.
But we already have performed
the symplectic reduction with respect to ${\cal G}_1$.
Therefore, we can push  $\om$ (\ref{16.2})
down on the factor space ${\cal P}''/{\cal G}_0$.
Since ${\cal G}_0$ acts on
${\cal N}'$ only,  it can be done by  fixing the dependence
of $\mu$ on the coordinates
 in the  Teichm\"{u}ller space
 ${\cal T}_{g,n}$.
According to (\ref{2b.2}) represent $\mu$ as
\beq{17a.2}
\mu=\sum_{s=1}^{l}\mu_s.
\eq
The Beltrami differential (\ref{17a.2}) defines the
tangent vector
$
{\bf t}=(t_1,\ldots,t_l),
$
to the Teichm\"{u}ller  space  ${\cal T}_{g,n}$ at the
fixed point of ${\cal T}_{g,n}$.

We specify the dependence of $\mu$ on the positions of
the marked points in the following  way.
Let ${\cal U}'_a\supset{\cal U}_a$ be two vicinities
of the marked point $x_a$
such that ${\cal U}'_a\cap{\cal U}'_b=\emptyset$ for $a\neq b$, and
 $\chi_a(z,\bz)$ is a smooth function
$$
\chi_a(z,\bz)=\left\{
\begin{array}{cl}
1,&\mbox{$z\in{\cal U}_a$ }\\
0,&\mbox{$z\in\Si_{g,n}\setminus {\cal U}'_a.$}
\end{array}
\right.
$$
Introduce times related to the positions of the
marked points $t_a=x_a-x_a^0$. Then
\beq{17.2}
\mu_a=t_a\bp n_a(z,\bz),~~n_a(z,\bz)=(1+c_a(z-x_a^0))\chi_a(z,\bz).
\eq
The action of ${\cal G}_0$ on the phase space ${\cal P}''$
reduces the infinite-dimensional component ${\cal N}'$ to
${\cal T}_{g,n}$.
After the reduction we come to the bundle with base
${\cal T}_{g,n}$.
The symplectic form (\ref{16.2}) is transformed as follows
\beq{19.2}
\om=\om_0({\bf v},{\bf u},{\bf p},{\bf t})-
\frac{1}{\ka}\sum_{s=1}^{l}\de H_s({\bf v},
{\bf u},{\bf p},{\bf t})\de t_s,~
H_s=\int_\Si <L,L>\p_s\mu
\eq
where $\om_0$ is defined by (\ref{15.2}).

In fact, we still have a redundant discrete  symmetry, since $\om$
 is invariant under the
mapping class group  $\pi_0({\cal G}_0)$. Eventually, we come to
the moduli space
 ${\cal M}_{g,n}={\cal T}_{g,n}/\pi_0({\cal G}_0)$.

The extended phase space ${\cal P}$ is the result of
the symplectic reduction
with respect to the ${\cal G}_1$ action and subsequent
factorization under
the ${\cal G}_0$ action.
We can write symbolically
${\cal P}=({\cal P}''//{\cal G}_1)/{\cal G}_0.$
It is endowed with the symplectic form (\ref{19.2}).
\bigskip

{\bf 5. The hierarchies of the isomonodromic
deformations (HID)}.
The equations of motion (HID) can be extracted from
the symplectic form
(\ref{19.2}).
In terms of the local coordinates
they take the form
\beq{21.2}
\ka\p_s{\bf v}=\{H_s,{\bf v}\}_{\om_0},~~
\ka\p_s{\bf u}=\{H_s,{\bf u}\}_{\om_0},~~
\ka\p_s{\bf p}=\{H_s,{\bf p}\}_{\om_0}~~
\eq
The Poisson bracket $\{\cdot,\cdot\}_{\om_0}$ is the
inverse tensor to $\om_0$.
We also has the  Whitham hierarchy  accompanying (\ref{21.2})
\beq{22.2}
\p_sH_r-\p_rH_s+\{H_r,H_s\}_{\om_0}=0.
\eq
There exists the one form on ${\cal M}_{g,n}$ defining {\sl the tau
function} of the hierarchy of isomonodromic deformations
\beq{22a.2}
\de \log\tau=\de^{-1}\om_0-\f1{\ka}\sum H_sdt_s.
\eq

The following three statements are valid for the HID (\ref{21.2}):
\begin{predl}
There exists the consistent system of
linear equations
\beq{23.2}
(\ka\p+L)\Psi=0,
\eq
\beq{24.2}
(\p_s+M_s)\Psi=0,~~(s=1,\ldots,l=\dim{\cal M}_{g,n})
\eq
\beq{25a.2}
(\bp+\mu\p+{\bar L})\Psi=0
\eq
 where $M_s$ is a solution to the linear equation
\beq{25.2}
\p_{\bar w}M_s-[M_s,{\bar L}]=\p_s\bar{L}-\frac{1}{\ka}L\p_s\mu.
\eq
\end{predl}

\begin{predl}
The linear conditions (\ref{24.2})
 provide the isomonodromic deformations of the linear system
 (\ref{23.2}), (\ref{25a.2})
 with respect to change the "times" on ${\cal M}_{g,n}$.
\end{predl}
Therefore, the HID (\ref{21.2}) are the monodromy preserving
conditions for the linear system (\ref{23.2}),(\ref{25a.2}).

The presence of derivative with
respect to the spectral parameter $w\in\Si_{g,n}$ in
the linear equation (\ref{23.2}) is a
distinguish feature of the monodromy preserving equations.
It plagues the
application of the inverse scattering method to these types
of systems. Nevertheless, in our case
 we have in some sense the explicit form of solutions:
\begin{predl}[The projection method.]
The solution of the Cauchy problem of (\ref{21.2})
 for the initial data  ${\bf v}^0,{\bf u}^0,{\bf p}^0$
 at the time
 ${\bf t}={\bf t}^0$
 is defined in terms of the elements $L^0,{\bar L}^0$
 as the gauge
transform
\beq{26a.2}
{\bar L}({\bf t})=f^{-1}( L^0(\mu({\bf t})-\mu({\bf t}^0))+
({\bar L}^0))f+
f^{-1}(\bp+\mu({\bf t})\p)f,
\eq
\beq{26.2}
L({\bf t})=f^{-1}(\p+L^0)f,~~{\bf p}({\bf t})=f^{-1}({\bf p}^0)f,
\eq
 where $f$  is a smooth  $G$-valued functions on $\Si_{g,n}$
fixing the gauge.
\end{predl}

\section{Relations to the Hitchin systems and the KZB equations.}
\setcounter{equation}{0}

{\bf 1. Scaling limit.}
Consider the HID in the limit $\ka\rar 0$.
We will prove that
in this limit we come to
the Hitchin systems, which are living on the cotangent
bundles to the
moduli space of holomorphic $G$-bundles over $\Si_{g,n}$
\cite{H1}.
The  critical value $\ka=0$ looks singular (see (\ref{2.2}),
(\ref{19.2})).
To get around we rescale the times ${\bf t}=\ka{\bf t}^H$,
where  ${\bf t}^H$ are the "Hitchin times". Therefore,
$
\de\mu({\bf t})=\ka\sum_s\p_s\mu({\bf t}^0)\de t^H.
$
After this rescaling the forms  (\ref{2.2}),(\ref{19.2})
become regular in
the critical limit.
The rescaling procedure means that we blow up a vicinity
of the fixed point
corresponding to $(\mu=0)$, and the whole dynamic of
the Hitchin systems
is developed in this vicinity.
\footnote{We are grateful to A.Losev for elucidating this point.}
 Denote $\p_s\mu_o=\p_s\mu(\bft)|_{\bft=\bft^0}$ 
Then we have instead of (\ref{2.2})
\beq{2.3}
\om=\int_{\Si}<\de A,\de\bar{A}>+
2\pi i\sum_{a=1}^n\de<p_a,g_a^{-1}\de g_a>-
\sum_s\int_{\Si}<\de A, A>\p_s\mu({\bf t}^0)\de t^H.
\eq
If $\ka=0$ the connection $A$ behaves  as the one-form
$A\in {\cal A}^{(1,0)}(\Si_{g,n},{\rm Lie}(G))$ (see (\ref{8.2})).
It is so called the Higgs field . An important point is that
the Hamiltonians now become the times independent. The form
(\ref{2.3}) is
the starting point in the derivation of the Hitchin systems via
the symplectic
reduction \cite {H1,Ne}. Essentially, it is the  same  procedure
as described above. Namely, we obtain the same moment constraint
(\ref{14.2})
and the same gauge fixing (\ref{11.2}).
But now we are sitting in a fixed point $\mu({\bf t}^0)=0$
of the moduli space ${\cal M}_{g,n}$
and don't need the factorization under the action of the
diffeomorphisms.
This only difference between the solutions
$L$ and $\bar{L}$  in the Hitchin systems
and the hierarchies of isomonodromic deformations.

Propositions 2.1, 2.3 are valid for the Hitchin systems in a
slightly modified form.
\begin{predl}
 There exists the consistent system of
linear equations
\beq{3.3}
(\la+L)\Psi=0,~~\la\in\bf{C}
\eq
\beq{4.3}
(\p_s+M_s)\Psi=0,~~\p_s=\p_{t^H_s},~(s=1,\ldots,l=\dim{\cal M}_{g,n})
\eq
\beq{5.3}
(\bp+{\bar L})\Psi=0, ~\bp=\p_{\bz},
\eq
 where $M_s$ is a solution to the linear equation
\beq{6.3}
\bp M_s-[M_s,{\bar L}]=\p_s\bar{L}-L\p_s\mu({\bf t}^0).
\eq
\end{predl}
The parameter $\la$ in (\ref{3.3}) can be considered as the symbol
of $\p$ (compare with (\ref{23.2})).

When $L$ and $M$ can be find explicitly the simplified
form of (\ref{14.2})
allows to apply "the inverse scattering method" to find  solutions
of the Hitchin hierarchy as it was done for $\SL$ holomorphic
bundles over $\Si_{1,1}$ \cite{Kr2}, corresponding to the elliptic
Calogero system with spins. We present the alternative way 
to describe the solutions:
\begin{predl}[The projection method.]
$$
\bar{L}(t_s)=f^{-1}(L^0\p_s\mu_o(t_s-t_s^0)+\bar{L}^0)f+f^{-1}\bp f,
$$
$$
L(t_s)=f^{-1}L^0f,~~p_a(t_s)=f^{-1}(p_a^0)f
$$
\end{predl}
The degenerate version of these expressions was known
for a long time \cite{OP}.
\bigskip
\noindent

{\bf 2. About KZB.}
The Hitchin systems are the classical limit of the KZB equations
on the critical level \cite{Ne,I}. The later has the form of
the Schr\"{o}dinger
equations, which is the result of geometric quantization of the moduli
of flat $G$
bundles \cite{ADPW,H2}. The conformal blocks of the WZW theory on
$\Si_{g,n}$
with vertex operators in marked points are ground state wave functions
$$\hat {H}_sF=0,~~(s=1,\ldots,l=\dim{\cal M}_{g,n}).$$
The classical limit means that one replaces operators on
their symbols and
finite-dimensional representations in the vertex operators
by the corresponding
coadjoint orbits. The level $\ka$ plays the role of
the Planck constant, but in
contrast with the limit considerd before, we don't adjust the moduli
of complex structures.

Generically, for $\ka\neq 0$ the KZB equations can be written
in the form of the nonstationar  Schr\"{o}dinger equations
$$
(\ka\p_s+\hat {H}_s)F=0,~~(s=1,\ldots,l=\dim{\cal M}_{g,n}).
$$
The flatness of this connection (see \cite{H2})
is the quantum counterpart of the Whitham equations (\ref{22.2}).
The classical limit in the described above sense leads to the HID.
 Summarizing, we arrange these
quantum and classical systems in the diagram.  The vertical
arrows denote to the classical limit, while  the limit
$\ka\rar 0$ on the horizontal arrows includes also the rescaling of
the moduli of complex structures. The examples in the bottom of the diagram
will be considered in next sections.
$$
\def\normalbaselines{\baselineskip20pt
       \lineskip3pt    \lineskiplimit3pt}
\def\mapright#1{\smash{
        \mathop{\longrightarrow}\limits^{#1}}}
\def\mapdown#1{\Big\downarrow\rlap
       {$\vcenter{\hbox{$\scriptstyle#1$}}$}}
\begin{array}{ccc}
\left\{
\begin{array}{c}
\mbox{KZB eqs.},~(\ka,{\cal M}_{g,n},G)\\
(\ka\p_{t_a}+\hat{H}_a)F=0,\\
(a=1,\ldots,\dim{\cal M}_{g,n})
\end{array}
\right \}
  &\mapright{\ka\rar0,~{\bf t}=\ka{\bf t}^H} &
\left\{
\begin{array}{c}
\mbox{KZB eqs. on the critical level},~\\
({\cal M}_{g,n},G),~(\hat{H}_a)F=0,\\
(a=1,\ldots,\dim{\cal M}_{g,n})
\end{array}
\right\}
\\
\mapdown{\ka\rar 0}&    &\mapdown{\ka\rar 0} \\
\left\{
\begin{array}{c}
\mbox{Hierarchies of Isomonodromic} ~\\
\mbox{deformations on}~{\cal M}_{g,n}
\end{array}
\right \}
&\mapright{\ka=0,~{\bf t}=\ka{\bf t}^H} &
\left\{
\begin{array}{c}
\mbox{Hitchin systems} ~\\
       \\
 \end{array}
\right \}
\\
 & & \\
  & \mbox{\sl EXAMPLES}& \\
\left\{
\begin{array}{c}
\mbox{Schlesinger eqs.} \\
\mbox{Painlev\'{e} type eqs.} \\
\mbox{Elliptic Schlesinger eqs.}
\end{array}
\right\}
  &  \mapright{\ka\rar 0,~{\bf t}=\ka{\bf t}^H} &
\left\{
\begin{array}{c}
\mbox{Classical Gaudin eqs.}\\\
\mbox{Calogero eqs.}\\
\mbox{Elliptic Gaudin eqs.}
\end{array}
\right\}
\end{array}
$$

\section{Genus zero - Schlesinger's equation.}
\setcounter{equation}{0}

Consider ${\bf C}P^1$ with $n$ punctures $(x_1,\ldots,x_n|x_a\neq x_b)$.
The Beltrami
differential $\mu$ is related only to the positions of marked points
(\ref{17.2}). On ${\bf C}P^1$  the gauge
transform (\ref{11.2}) allows to choose $\bar L$ to be identically zero.
Let $A=f(L+\ka\p_w)f^{-1}$. Then
the moment equation takes the form
\beq{2.1}
(\bp +\p\mu)L=2\pi i\sum_{a=1}^n\de^2(x_a)p_a.
\eq
 It allows to find $L$
\beq{2.5}
L=\sum_{a=1}^n\frac{p_a}{w-x_a}.
\eq
On the symplectic quotient $\om$ (\ref{19.2})
takes the form
$$\om=\de\sum_{a=1}^n<p_ag_a^{-1}\de g_a>-
\frac{1}{\ka}\sum_{b=1}^n(\de H_{b,1}+\de H_{b,0})\de x_b.
$$
\beq{3.5}
H_{a,1}=\sum_{b\neq a}\frac{<p_a,p_b>}{x_a-x_b},~~H_{2,a}=c_a<p_a,p_a>.
\eq
$H_{1,a}$ are
precisely the Schlesinger's Hamiltonians.
Note, that we still have a gauge freedom with respect to the $G$ action.
The corresponding moment constraint means that the sum of residues of
$L$ vanishes:
\beq{VP}
\sum_{a=1}^np_a=0.
\eq
While $H_{2,a}$ are Casimirs and lead to trivial equations,
the equation of motion for $H_{1,a}$ are the Schlesinger equations
$$
\ka\p_bp_a=\frac{[p_a,p_b]}{x_a-x_b},~(a\neq b),~~
\ka\p_ap_a=-\sum_{b\neq a}\frac{[p_a,p_b]}{x_a-x_b}.
$$
As by product, we obtain by this procedure the corresponding
 linear problem
(\ref{23.2}),(\ref{24.2})
with
$L$ (\ref{2.5}) and
$$M_{a,1}=-\frac{p_a}{w-x_a}$$
 as a solution to (\ref{25.2}).
The tau-function for the Schlesinger equations has the form \cite{JMU}
$$\de\log\tau=\sum_{c\neq b}<p_b,p_c>\de\log (x_c-x_b).$$

\section{Genus one - elliptic Schlesinger, Painleve VI...}
\setcounter{equation}{0}

{\bf 1. Deformations of elliptic curves.}
In addition to the moduli coming from the positions of the marked points
there is an elliptic module $\tau,~Im\tau>0$ on $\Si_{1,n}$.
As in (\ref{17a.2}),(\ref{17.2}) we take the Beltrami differential 
in the
form \mbox{$
\mu=\sum_{a=1}^n\mu_a+\mu_\tau,~~(\mu_s=t_s\bp n_s),
$}
where $n_a(z,\bz)$ is the same as in (\ref{17.2}) and
\beq{1.6}
n_\tau=(\bz-z)(1-\sum_{a=1}^n\chi_a(z,\bz)).
\eq
Then
\beq{2b.6}
 \mu_\tau=
\ti{\mu}_\tau(1-\sum_{a=1}^n\chi_a(z,\bz)),~~
(\ti{\mu}_\tau=\frac{t_\tau}{\tau-\tau_0},~t_\tau=\tau-\tau_0)
\eq
Here $\tau_0$ defines the reference comlex structure on the curve
$$
T^2_{0}=\{0<x\leq 1,~0<y\leq 1, ~z=x+\tau_0y,~\bz=x+\bar{\tau}_0y\}.
$$

\bigskip
{\bf 2. Flat bundles on a family of elliptic curves.}
Note first, that $\bA$ can be considered as a connection of holomorphic
$G$ bundle ${\cal E}$ over $T^2_{\tau}$.
For stable bundles $\bA$ can be gauge transformed by
(\ref{11.2}) to the Cartan $z$-independent form\\
$\bA=f(\bar{L}+\bp+\mu\p)f^{-1}$,
$ \bar{L}\in{\cal H}$- Cartan subalgebra of Lie$(G)$.
Therefore, stable bundle ${\cal E}$ is decoposed into the direct
sum of line bundles
${\cal E}=\oplus_{k=1}^r {\cal L}_k, ~~r=\mbox{rank}G$.
The set of gauge equivalent connections represented by $\{\bar{L}\}$
can be identified with the $r$ power of the Jacobian of $T^2_{\tau}$,
factorized by the action of the Weyl group $W$ of $G$.
Put
\beq{4.6}
\bar{L}=2\pi i \frac{1-\ti{\mu}_\tau}{\rho}{\bf u},
 ~~{\bf u}\in{\cal H},
~~(\rho=\tau_0-\bar{\tau}_0).
\eq
The moment constraints (\ref{13.2}) leading to the flatness
condition take the form
\beq{3.6}
\p_{\bar{w}}L+[\bar{L},L]=2\pi i\sum_{a=1}^n\de^2(x_a,)p_a.
\eq
Let $R=\{\al\}$ be the root system of of Lie$(G)={\cal G}$ and
${\cal G}={\cal H}\oplus_{\al\in R}{\cal G}_\al$
be the root decomposition.
Impose the vanishing of the residues in (\ref{3.6})
\beq{11.6}
\sum_{a=1}^n(p_a)_{\cal H}=0,
\eq
where $p_a|_{\cal H}$ is the Cartan component and we have
identified ${\cal G}$ with its dual space ${\cal G}^*$.

We will parametrized the set of its solutions by two elements
$\bfv,{\bf u}\in {\cal H}$.
Define the solutions $L$ to the moment equation (\ref{3.6}),
which is double periodic on the deformed curve  $T^2_\tau$. Let
$E_1(w)$ be the Eisenstein function
$$
E_1(z|\tau)=\p_z\log\te(z|\tau),
$$
where
$$
\te(z|\tau)=q^{\frac
{1}{8}}\sum_{n\in {\bf Z}}(-1)^ne^{\pi i(n(n+1)\tau+2nz)}.
$$
It is connected with the Weirstrass zeta-function as
$$
\ze(z|\tau)=E_1(z|\tau)+2\eta_1(\tau)z,~~(\eta_1(\tau)=
\ze(\frac{1}{2})).
$$
Another function we need is
$$
\phi(u,z)=\frac{\te(u+z)\te'(0)}{\te(u)\te(z)}=
\exp(-2\eta_1uz)
\frac
{\si(u+z)}{\si(y)\si(z)},
$$
where $\si(z)$ is the Weierstrass sigma function.
\begin{lem}
The solutions of the moment constraint equation have the form
\beq{5.6}
L=P+X,~~P\in{\cal H},~~X=\sum_{\al\in R}X_{\al}.\eq
\beq{6.6}
P=2\pi i(\frac{\bfv}{1-\ti{\mu}_\tau}-\ka\frac{{\bf u}}{\rho}+
\sum_{a=1}^n(p_a)_{\cal H}E_1(w-x_a)),
\eq
\beq{7.6}
X_{\al}=frac{2\pi i}{1-\ti{\mu}_\tau}
\sum_{a=1}^n(p_a)_{\al}\exp 2\pi i\{
\frac{(w-x_a)-(\bar{w}-\bar{x}_a)}{\tau-\bar{\tau}_0}
\al(u)\}
\phi(\al(u),w-x_a).
\eq
\end{lem}
\bigskip
{\bf 2. Symmetries.}
The remnant gauge transforms preserve the chosen Cartan subalebra
${\cal H}\subset G$. These transformations are generated by the Weyl
 subgroup $W$ of $G$ and  elements
\mbox{$f(w,\bar{w})\in {\rm Map}(T^2_{\tau},{\rm Cartan}(G))$}.
Let $\Pi$ be the system of simple roots,
$R^{\vee}=\{\al^{\vee}=\frac{2\al}{(\al|\al)}\},$ is the
dual root system,
and ${\bf m}=\sum_{\al\in \Pi} m_{\al}\al^{\vee}$
 be the element from the dual root
lattice ${\bf Z}R^{\vee}$. Then the Cartan valued harmonics
\beq{9.6}
f_{{\bf m},{\bf n}}=\exp 2\pi i(
{\bf m}\frac{w-\bar{w}}{\tau-\tau_0}+
{\bf n}\frac{\tau\bar{w}
-\bar{\tau}_0w}{\tau-\tau_0}),~~
({\bf m} ,{\bf n}\in R^{\vee})
\eq
generate the basis in the space of Cartan gauge transformations.
In terms of the variables $\bfv$ and ${\bf u}$
 they act as
\beq{10.6}
{\bf u}\rar {\bf u} +{\bf m}-{\bf n}\tau,~~ \bfv\rar
\bfv-\ka{\bf n},~~
(p_a)_{\al}\rar\vf(m_\al,n_\al)(p_a)_{\al}.
\eq
Here
$
\vf(m_\al,n_\al)=
\exp \frac{4\pi i}{\rho}[(m_\al-n_\al\bar{\tau}_0x_a^0)-
(m_\al-n_\al\tau_0\bar{x}_a^0)].$

The whole discrete gauge symmetry is the semidirect product
$\hat{W}$ of
 the Weyl group
$W$ and the lattice ${\bf Z}R^{\vee}\oplus\tau{\bf Z}R^{\vee}$.
It is the
Bernstein-Schvartsmann complex crystallographic group.
The factor space ${\cal H}/\hat{W}$ is  the genuin space for
the "coordinates" ${\bf u}$.

 According with (\ref{8.2}) the transformations (\ref{9.6})
 act also on
$p_a\in{\cal O}_a$. This action leads to the symplectic quotient
${\cal O}_a//H$ and generates the moment equation (\ref{11.6}).

The modular group ${\rm PSL}_2({\bf Z})$ is a subgroup of mapping
class group
for the Teichm\"{u}ller space ${\cal T}_{1,n}$.
Its action on $\tau$ is the M\"{o}bius transform.
We summarise the action of symmetries on the dynamical
variables:
\bigskip
\begin{center}
\begin{tabular}{|c|c|c|c|}\hline
             &W=\{s\}&${\bf Z}R^{\vee}\oplus\tau{\bf Z}R^{\vee}$
             &${\rm PSL}_2({\bf Z})$
\\ \hline \hline
$\bfv$ &$s\bfv$& $\bfv-\ka{\bf n} $ & $\bfv(c\tau+d)-\ka c{\bf u} $
\\ \hline
${\bf u}$ &$s{\bf u}$  & ${\bf u}+{\bf m}-{\bf n}\tau$
&${\bf u}(c\tau+d)^{-1}$
\\ \hline
$(p_a)_{\cal H}$ & $s(p_a)_{\cal H}$   &   $(p_a)_{\cal H}$
& $(p_a)_{\cal H}$
\\   \hline
$(p_a)_{\al}$ & $(p_a)_{s\al}$   &   $\vf(m_\al,n_\al)(p_a)_{\al}$
& $(p_a)_\al$
\\   \hline
$\tau$   &$\tau$  &$\tau$                   &$\frac{a\tau+b}{c\tau+d}$
\\ \hline
$x_a$     &  $x_a$  &$x_a$        &    $\frac {x_a}{c\tau+d}$
\\ \hline
\end{tabular}
\end{center}

\bigskip
{\bf 3. Symplectic form.}
The set $(\bfv,{\bf u}\in {\cal H},{\bf p}=(p_1,\ldots,p_n))$
of dynamical variables along with the times
${\bf t}=(t_\tau,t_1,\ldots,t_n)$
describe the local coordinates in the bundle ${\cal P}$. According with
the general prescription, we can define the hamiltonian system
 on this set.

The main statement, formulated in terms of the theta-functions and the
Eisenstein functions
\beq{A.2}
E_2(z|\tau)=-\p_zE_1(z|\tau)=
\p_z^2\log\te(z|\tau)=\wp(z|\tau)+2\eta_1(\tau).
\eq
It takes the form
\begin{predl}
The symplectic form $\om$ (\ref{19.2}) on ${\cal P}$ is
\beq{16.6}
\frac{1}{4\pi^2}\om=(\de\bfv,\de{\bf u})+
\sum_{a=1}^n\de<p_a,g_a^{-1}\de g_a>
-\frac{1}{\ka}(\sum_{a=1}^n\de H_{2,a}+\de H_{1,a})\de t_a-
\frac{1}{\ka}\de H_{\tau}\de\tau,
\eq
with the Hamiltonians
$$H_{2,a}=c_a<p_a,p_a>;$$
$$
H_{1,a}=
=2(\frac{\bfv}{1-\ti{\mu}_\tau}-\ka\frac{{\bf u}}{\rho},p_a|_{\cal H})+
\sum_{b\neq a}(p_a|_{\cal H},p_b|_{\cal H})E_1(x_a-x_b)+
$$
$$
\sum_{b\neq a}\sum_\al(p_a|_{\al},p_b|_{-\al})
\frac
{\te(-\al({\bf u})+x_a-x_b)\te'(0)}
{\te(\al({\bf u}))\te(x_a-x_b)};
$$
$$
H_{\tau}=
$$
$$\frac{(\bfv,\bfv)}{2}+\{\sum_{a=1}^n
\sum_{\al}(p_a|_{\al},p_a|_{-\al})E_2(\al(\bfu))+
\sum_{a\neq b}^n(p_a|_{\cal H},p_b|_{\cal H})(E_2(x_a-x_b)-
E_1^2(x_a-x_b))+
$$
$$
\sum_{a\neq b}\sum_\al(p_a|_{\al},p_b|_{-\al})
\frac
{\te(-\al(\bfu)+x_a-x_b)\te'(0)}
{\te(\al(\bfu))\te(x_a-x_b)}
(E_1(\al(\bfu)-E_1(x_b-x_a+\al(\bfu))-E_1(x_b-x_a)).
$$
\end{predl}
\bigskip
{\bf Example 1.}
Consider ${\rm SL}(2,{\bf C})$ bundles over the family of $\Si_{1,1}$.
Then (\ref{4.6}) takes the form
\beq{BL}
\bar{L}=2\pi i \frac{1-\ti{\mu}_\tau}{\rho}\di (u,-u).
\eq
In this case the position of the maked point is no long the module and we
put $x_1=0$. Since $\dim{\cal O}=2$ the orbit degrees of freedom
can be gauged away by the hamiltonian action of the diagonal group.
We assume that
$p=\nu[(1,1)^T\otimes(1,1)-Id]$.
Then we have from(\ref{5.6}),(\ref{6.6}),(\ref{7.6})
\beq{17.6}
L=2\pi i \mat{\frac{v}{1-\ti{\mu}_\tau}-\ka\frac{u}{\rho}}
{x(u,w,\bar{w})}{x(-u,w,\bar{w})}
{-\frac{v}{1-\ti{\mu}_\tau}+\ka\frac{u}{\rho}}.
\eq
$$
x(u,w,\bar{w})=\frac{\nu}{1-\ti{\mu}_\tau}\exp 4\pi i\{(w-\bar{w})u
\frac{1-\ti{\mu}_{\tau}}{\rho}\}\phi(2u,w).
$$

The symplectic form (\ref{16.6})
$$
\frac{1}{4\pi^2}\om=(\de v,\de u)-\frac{1}{\ka}\de H_{\tau}\de\tau,
$$
and
$$
H_{\tau}=v^2+U(u|\tau),~U(u|\tau)=-\nu^2E_2(2u|\tau).
$$
It leads to the equation of motion
\beq{20.6}
\frac{\p^2 u}{\p\tau^2}=\frac{2\nu^2}{\ka^2}\frac{\p}{\p u}E_2(2u|\tau).
\eq
In fact, due to (\ref{A.2}) we can use $\wp(2u|\tau)$ instead of 
$E_2(2u|\tau)$
and after rescaling the coupling constant come to (\ref{I.3}) for special
values of constants as in (\ref{I.8}).
 The equation  (\ref{20.6}) is the isomonodromic deformation conditions 
for the linear system
(\ref{23.2}),(\ref{25a.2}) with $L$ (\ref{17.6}) and $\bar{L}$ (\ref{BL}).
The lax pair is given by $L$ (\ref{17.6}) and $M_\tau$
$$
M_\tau=\mat{0}{y(u,w,\bar{w})}{y(-u,w,\bar{w})}{0},
$$
where $y(u,w,\bar{w})$ is defined as the convolution integral on $T^2_\tau$
$$y(u,w,\bar{w})=-\f1{\ka}x\ast x(u,w,\bar{w}).$$
The projection method determines  solutions of (\ref{20.6}) as a result of
diagonalization of $L$ (\ref{17.6}) by the gauge transform on the deformed
curve $T^2_\tau$.
On the critical level $(\ka=0)$ we come to the two-body elliptic Calogero
 system.
\bigskip

{\bf Example 2.}
For  flat $G$ bundles over $\Si_{1,1}$ we obtain Painlev\'{e}
 type equations,
related to arbitrary root system. They are described by the system
of differential equations for the
$\bfu=(u_1,\ldots,u_r),~(r=$rank$G)$ variables.
In addition there are the orbit variables $p\in{\cal O}(G)$
satisfying the Euler top equations. For $\SL$
bundles  the most degenerate orbits ${\cal O}=T^*{\bf C}P^{N-1}$ 
has dimension
$2N-2$. These variables are gauge away by the diagonal gauge transforms as in
the previous example. On the critical level this Painlev\'{e} type
system degenerates into $N$-body elliptic Calogero
system. For generic orbits  we obtain the generalized Calogero-Euler systems.

\small{

}
\end{document}